%% file: ms.tex
\documentclass[twocolumn]{autart_mod_arXiv}    

\usepackage{hyperref, bookmark}
\hypersetup{colorlinks=false,pdfborder=0 0 0}

\usepackage{amssymb}
\usepackage{graphicx,setspace}
\usepackage{amsmath,amssymb}
\usepackage{color}
\usepackage[english]{babel}
\usepackage{mathrsfs,dsfont}
\usepackage{array}
\usepackage{tikz}
\usepackage{adjustbox}
\usepackage{multirow}
\usetikzlibrary{arrows, automata, decorations.pathmorphing, backgrounds, positioning, fit, petri, matrix, patterns,decorations.pathreplacing}

\usepackage{algorithm}
\usepackage{etoolbox}
\let\classAND\AND
\let\AND\relax
\usepackage{algorithmic}

\let\AND\classAND
\AtBeginEnvironment{algorithmic}{\let\AND\algoAND}

\usepackage[detect-none]{siunitx}
\usepackage{lineno}

\usepackage{mathrsfs}
\usepackage[mathscr]{euscript}

\def\QED{~\rule[-1pt]{5pt}{5pt}}

\let\VEC  \boldsymbol
\let\ALP  \mathcal

\newcommand{\transpose}{\mathsf{T}}

\newcommand{\beq}[1]{\begin{linenomath}\begin{align} #1 \end{align}\end{linenomath}}
\newcommand{\beqq}[1]{\begin{linenomath}\begin{align*} #1 \end{align*}\end{linenomath}}
\renewcommand{\Re}{\mathbb{R}}

\begin{document}

\begin{frontmatter}

%

\title{An Algorithm to Warm Start Perturbed (WASP) Constrained Dynamic Programs} 

\author[1]{Abhishek Gupta }\ead{gupta.706@osu.edu},    
\author[2]{Shreshta Rajakumar Deshpande}\ead{rajakumardeshpande.1@osu.edu},               
\author[2]{Marcello Canova}\ead{canova.1@osu.edu}  

\address[1]{Department of Electrical and Computer Engineering, The Ohio State University, Columbus, OH, USA - 43210}  
\address[2]{Department of Mechanical and Aerospace Engineering, The Ohio State University, Columbus, OH, USA - 43210}


\begin{abstract}                          
Receding horizon optimal control problems compute the solution at each time step to operate the system on a near-optimal path. However, in many practical cases, the boundary conditions, such as external inputs, constraint equations, or the objective function, vary only marginally from one time step to the next. In this case, recomputing the optimal solution at each time represents a significant burden for real-time applications. This paper proposes a novel algorithm to approximately solve a perturbed constrained dynamic program that significantly improves the computational burden when the objective function and the constraints are perturbed slightly. The method hinges on determining closed-form expressions for first-order perturbations in the optimal strategy and the Lagrange multipliers of the perturbed constrained dynamic programming problem are obtained. This information can be used to initialize any algorithm (such as the method of Lagrange multipliers, or the augmented Lagrangian method) to solve the perturbed dynamic programming problem with minimal computational resources.  
\end{abstract}

\end{frontmatter}

\section{Introduction}
Many real-time applications require solving receding horizon constrained dynamic optimization problems repeatedly to operate the system on a near-optimal path. As an example, a typical problem in Connected and Autonomous Vehicles (CAVs) is ``Eco-Driving", where the vehicle velocity is optimized dynamically to determine the most fuel-efficient trajectory based upon the route characteristics (e.g., speed limits, grade) and the presence of other vehicles. This problem has been approached using a variety of control methods including Pontryagin's Minimum Principle (PMP) in \cite{ozatay2017velocity,uebel2017optimal}, Dynamic Programming (DP) in \cite{olin2019reducing,deshpande2020integrated,heppeler2014fuel,kuriyama2010theoretical}, Machine Learning methods in \cite{liu2017reinforcement,liu2019adaptive}, and Model Predictive Control (MPC) \cite{homchaudhuri2015fuel,homchaudhuri2016fast,sun2014velocity,yamaguchi2012model}. During the vehicle operation, it is often the case that parameters in the state dynamics (such as vehicle velocity and battery state-of-charge) or objective function (such as the fuel consumption), or even constraint equations (such as the route speed limits and actuator operating limits), may change as a result of variations in the environmental conditions or external inputs (such as the presence of a lead vehicle).

Another real-time application is the problem of scheduling the charging process of electric vehicles (EVs) using renewable energy. An emerging approach to solve this problem is using MPC, that takes into account prediction for the renewable energy generation \cite{di2014electric,chen2014distributional,lee2018large}. For illustration, consider the following problem setup. The state of the system is the battery state-of-charge of the EVs and the action is the amount of energy allocated to each EV in each time step. The cost function of the system is to charge all the EVs connected to the system with minimal cost of energy acquired from non-renewable sources. The constraints are typically total energy constraints and the charging rate constraints of the EVs. The state dynamics is typically linear in the state and the action. As time progresses, the renewable forecasts are updated, new EVs arrive, charged EVs depart, and the cost of electricity from non-renewable sources changes (which typically reflects the real-time market). The MPC problem is changed accordingly to reflect the new market conditions. Many of the changes illustrated in these examples can be captured by applying small perturbations to the objective function or constraints, without fundamentally altering the mathematical formulation of the Receding Horizon Optimal Control Problem (RHOCP).

In general, two approaches are adopted to compute the optimal solution to a constrained dynamic optimization problem. The most common method is to determine an optimal ``open-loop'' policy, in which the optimal actions are computed along with the optimal system trajectory, resulting into a weakly time-consistent optimal policy. If the system deviates from the optimal trajectory, then the entire optimal solution needs to be recomputed. The key mathematical tool to enable this class of algorithms is the maximum principle, which is essentially a gradient descent algorithm specialized to discrete-time dynamic systems; see the discussion in Section 2.6 of \cite{bertsekas2016nonlinear}.

The other approach to solve for the optimal solution is to use a DP algorithm, in which the optimal solution is computed backward in time \cite{bertsekas2012dynamic,bellman2015applied,larson1978principles}. This computationally heavy algorithm yields the optimal ``closed-loop'' policy, in which the optimal actions are computed as functions of the state, and the optimal cost-to-go functions using optimal closed-loop policies are kept track of using ``value functions''. This yields a strongly time-consistent optimal policy, that is, even if the system deviates from the optimal trajectory, one can just use the optimal policy and the current state to recompute the optimal actions for the current as well as the future time steps. The use of DP thus yields a policy that is robust to external disturbances.

However, in many practical cases, the number of states in the plant dynamics discourages from applying DP, due to the computational complexity. In the example of Eco-Driving for CAVs, the issue arises when the vehicle operates in ``car-following" mode. Since the controlled vehicle must maintain a safe distance with respect to the leading vehicle, the relative position or time must added to the vehicle velocity as an additional state. 
The issue becomes more complex by the fact that constraint equations could be included or removed at each time step (for instance, a car from the side lane enters the lane of the CAV) or changes in external conditions may impact the state dynamics and objective function (for instance, auxiliary loads may increase the load on the powertrain and, ultimately, the fuel consumption). The challenge of handling perturbations in the parameters or constraint equations in a systematic way warrants the need to develop methods that can use real-time information to update the RHOCP formulation and execute the DP algorithm under in a computationally efficient manner.


It is worth noting that, in many practical applications of optimal control, the external variables change gradually, and potentially at much slower rate than the time scale considered in the computation horizon. This is equivalent of considering a \textit{perturbed RHOCP}, in which the cost functions and the constraints only change slightly from the previous solution. 

This paper proposes a novel approach to speed up the computation of Dynamic Programming for solving a perturbed RHOCP. The approach, a Warm Start Programming (WASP) utilizes a previously computed solution as the starting point for the next iteration of the RHOCP formulated after the parameters, constraints, or external variables have changed. Simple matrix operations are used to determine the first-order perturbations in the optimal solutions, including the perturbations in optimal strategy and the Lagrange multipliers of the RHOCP, at each time step. Exact expressions for the first order perturbations in the optimal solution are derived, and the efficacy of the method is illustrated for two dynamic optimization problems.


\subsection{Prior Work}
The foundations of dynamic programming can be traced to the late 1940s in the confidential work of Isaacs \cite{isaacs1999differential} in the context of differential games. Rediscovered by Bellman in \cite{bellman1952theory}, this theory has been further developed for finite and infinite horizon optimization problems \cite{bertsekas2012dynamic}, Markov decision problems \cite{altman1999constrained}, and more recently, over abstract spaces \cite{bertsekas2018abstract,yuksel2020universal}. The theory has found numerous applications in control theory \cite{kumar2015stochastic}, operations research problems \cite{si2004handbook}, automotive systems \cite{borrelli2017predictive}, economics \cite{rust1996numerical}, among many others.

Perturbations of static optimization problem is a very well-studied problem. The continuity properties of the optimal solution and optimal value as functions of the parameters was studied in \cite{dantzig1967continuity}. For linear programs, \cite{bohm1975continuity,bereanu1976continuity} showed that under certain regularity conditions, the optimal solution and optimal value are continuous in the parameters. For stochastic linear programs with recourse, \cite{walkup1969stochastic} establishes the lower semicontinuity of the value function. The continuity result was later established for quadratic optimization problem with linear constraints in \cite{tam1999continuity,lee2006continuity}.

For parametrized convex programs, \cite{hogan1973continuity} establishes continuity properties of the optimal value as a function of certain specific perturbations. Reference \cite{li2014holder} showed H\"older continuity of the optimal solution with respect to the parameters in strongly convex programs with parametrized constraint sets satisfying certain properties. For uniquely solvable linear complementarity problems, \cite{mangasarian1987lipschitz} characterized a Lipschitz constant for the continuity of the unique optimal solution (again, as a function of the parameters). For optimization over Banach spaces, \cite{shapiro1992perturbation} identified a set of sufficient conditions on the optimization functionals over Banach spaces such that the optimal solutions are continuous in perturbations; see also the survey paper \cite{bonnans1998optimization} and the recent text \cite{bonnans2013perturbation} on related topics.

Simultaneous to the above development on static optimization problems, limited literature has studied the continuity of closed-loop policies in finite and infinite horizon dynamic programs. In particular, the continuity of both the value function and the optimal closed-loop policy in an infinite horizon dynamic program was established in \cite{jordan1977continuity}.

Building upon these results, this paper shows that if certain continuity and differentiability assumptions hold for the optimal closed-loop policies and the value functions, and if the dynamic program is perturbed slightly, then one can compute an approximately optimal solution by simple matrix manipulations. This can then be fed to a DP solver as the initial guess so that the overall RHOCP problem can be solved quickly and effectively in real-time applications. The approach developed here requires slightly more memory to store certain derivatives and Lagrange multipliers.


%

\subsection{Notation}
Let $f:\Re^n\rightarrow\Re$ be a differentiable function. Then, $\nabla_x f$ denotes the derivative of the function $f$, and is given by
\beqq{\nabla_x f(x) = \begin{bmatrix}\nabla_{x_1}f(x)\\ \vdots \\ \nabla_{x_n}f(x)\end{bmatrix}.}
Similarly, if $f:\Re^n\rightarrow\Re^m$ is a differentiable function such that $f(x) = [f_1(x),\ldots,f_m(x)]^\transpose$, then its derivative is 
\beqq{\nabla_x f(x) = \begin{bmatrix}\nabla_xf_1(x)|\cdots | \nabla_x f_m(x)\end{bmatrix}.}
Let $f:\Re^n\rightarrow \Re^m$ and $g:\Re^m\rightarrow\Re^p$ be differentiable functions. The chain rule for derivative of the composite function $g\circ f:\Re^n\rightarrow \Re^p$ is given as
\beq{\nabla_x g\circ f|_{x_0} = \nabla_x f|_{x_0} \nabla_f g|_{f(x_0)}. \label{eqn:chainrule}}

\subsection{Outline of the Paper}
The problem formulation is introduced in the next section. In Section \ref{sec:static}, simple perturbations to a quadratic optimization problem are first considered. The insights and results are then used to determine the first order perturbations in the optimal solution and corresponding Lagrange multipliers of a general nonlinear optimization problem. In Section \ref{sec:wasp}, this approximation scheme is bootstrapped for a dynamic optimization problem to yield the WASP algorithm. Here, two numerical examples are constructed to demonstrate results from the proposed WASP algorithm. Some simplifications to the algorithm are discussed in Section \ref{sec:discuss}, and conclusions are finally drawn in Section \ref{sec:conclusion}.

\section{Problem Formulation}\label{sec:problem}
\input{2_Section_Problem_Formulation}

\section{Warm Start Static Optimization}\label{sec:static}
\input{3_Section_Constraint_Perturbation_Warm_Start}

\section{Warm Start Dynamic Programming (WASP)}\label{sec:wasp}
\input{4_Section_WASP}

\section{Discussions}\label{sec:discuss}
In Theorem \ref{thm:dpperturb}, observe that the first order perturbation $d_t(x)$ needs to be computed at every $x$, which may not be difficult if all the expressions required in \eqref{eqn:dpperturb} are available in the closed form. However, this may not be the case in many practical applications. Thus, computing the first order perturbation in the optimal strategy may become as computationally demanding as computing the DP solution from scratch. This computation can be simplified as follows.

If the first order perturbation is not significant (as is the case when $\epsilon$ is very small), then it is reasonable to assume that the state trajectories under $\bar{\VEC \gamma}^*$ and $\VEC \gamma^*$ are close to each other. Let $x^*_t$ denote the trajectory under the influence of $\VEC\gamma^*$. Then, $d_t(x) = d_t:= B_t(x^*_t)w_t(x^*_t)$ for all $x\in\ALP X$. Similarly, when taking the derivative of $\tilde V_t\circ f_{t-1}$ with respect to $u$, certain derivative terms can be ignored if the corresponding terms do not have large curvature. Such simplifications are likely going to be problem dependent, so such details are left to the practitioners.

Some limitations of the WASP algorithm and possible remedies are summarized below:
\begin{enumerate}
	\item After discretizing the state space, the value functions are stored in a tabular format. The first and second derivatives are computed using the finite difference method. This can lead to propagation of inaccuracies in the value function in dynamic programs with large horizon $T$.
	
	One way to alleviate this issue is to use smooth function approximators like polynomials, neural networks, functions in a reproducing kernel Hilbert space, or support vector machines so that the value function is smooth and the computation of the derivatives is simplified.
	\item As previously noted, the WASP algorithm (in which $d_t(x)$ is computed using \eqref{eqn:dteq}) does not perform well when some constraints are active in certain regions of the state space and inactive in others.
	
	Using \eqref{eqn:dtineq} to compute $d_t(x)$ has demonstrated better approximations. However, this requires solving a quadratic optimization problem, whose computational complexity is slightly higher than the matrix multiplication used in \eqref{eqn:dteq}.

	\item These algorithms might require slightly larger memory than that required for running the usual DP algorithms (for storing Lagrange multipliers and derivatives of the cost and value functions).
\end{enumerate}

\section{Conclusion}\label{sec:conclusion}
In this paper, the perturbations of constrained dynamic programs were discussed. Expressions for the first order perturbations in the optimal strategy were derived, which depend on the optimal strategy and Lagrange multipliers for the original unperturbed problem. This framework yields fast model predictive control and approximate dynamic programming algorithms for real-time implementation with limited computational power.

As future work, this method can be extended to situations where some of the control strategies or state transition functions are in the form of look-up table or maps, and for cases where some of the actions take discrete values. More recently, learning-based control uses function approximators like neural networks to store the value function and the policy functions. It would also be interesting to explore the development of fast dynamic programming algorithms for such situations.

\begin{ack}                               
The work was supported in part by the United States Department of Energy, Advanced Research Projects Agency – Energy (award number DE-AR0000794), and in part by NSF ECCS Grant 1610615, which the authors gratefully acknowledge.  
\end{ack}

\bibliographystyle{plain}        
\bibliography{app,econ,math,opt,DP}

\appendix

\section{Proof of Lemma \ref{lem:quadapprox1}} \label{app:quadapprox1}

The inverse of a matrix $H$ perturbed by $\epsilon \tilde H$ can be approximated by \cite[Section 5.1.10.2]{zwillinger2002crc}:
\beq{(H+ \epsilon\tilde H)^{-1} = H^{-1} + \bar{H} + o(\epsilon), \label{eqn:approxpertmat}}
where $\bar{H} = -\epsilon H^{-1}\tilde H H^{-1}$. To see this, note that
\beqq{(H^{-1} + \bar{H} + o(\epsilon) )(H + \epsilon \tilde H) = I,}
where $I$ is the identity matrix of appropriate dimension. By comparing the first-order terms on both sides of this equation, we get
\beqq{\bar{H} = -H^{-1} \epsilon \tilde H H^{-1}.}
On substituting this result in \eqref{eqn:approxpertmat}, we arrive at \eqref{eqn:heps}. Now we use it to evaluate
\beqq{\Big(A(H&+\epsilon\tilde H)^{-1}A^\transpose\Big)^{-1}\\
	&= \Big(A(H^{-1} - H^{-1}\epsilon\tilde H H^{-1} + o(\epsilon))A^\transpose\Big)^{-1}\\
	&= (AH^{-1}A^\transpose - \epsilon AH^{-1}\tilde H H^{-1}A^\transpose + o(\epsilon) )^{-1}.}
Using the definitions of $M = (AH^{-1}A^\transpose)^{-1}$ and $\tilde M = AH^{-1}\tilde H H^{-1}A^\transpose$ and the result in \eqref{eqn:heps}, we get
\beq{\Big(A(H&+\epsilon\tilde H)^{-1}A^\transpose\Big)^{-1} = M + \epsilon M\tilde M M + o(\epsilon).}
In order to derive the expression for $\tilde \kappa^* - \kappa^*$, let us first note the following fact: If $C_1,C_2,\tilde C_1,\tilde C_2$ are matrices of appropriate dimensions, then 
\beq{(C_1 & +\epsilon \tilde C_1+o(\epsilon))(C_2+\epsilon \tilde C_2+o(\epsilon)) - C_1 C_2 \nonumber\\
	& = \epsilon(C_1\tilde C_2+\tilde C_1 C_2)+o(\epsilon).\label{eqn:c1c2}}
From the definition of $\kappa^*$ in \eqref{eqn:quadprogkappa}:
\beqq{\kappa^* = -(AH^{-1}A^\transpose)^{-1}(AH^{-1} e+b).}
Using the expressions in \eqref{eqn:c1c2}, we have
\beqq{& (H+\epsilon \tilde H)^{-1} (e+\epsilon \tilde e) - H^{-1} e \\
	& = (H^{-1}-\epsilon H^{-1}\tilde H H^{-1}+o(\epsilon))  (e+\epsilon \tilde e) - H^{-1} e\\
	& = \epsilon(H^{-1} \tilde e-H^{-1} \tilde H H^{-1} e)+o(\epsilon).}
This immediately yields
\beqq{& \Big(A(H+\epsilon \tilde H)^{-1} (e+\epsilon \tilde e)+b+\epsilon\tilde b\Big) - (AH^{-1} e+b) \\
	& = \epsilon\tilde m +o(\epsilon).}
Collecting the above expressions with \eqref{eqn:ahaeps} and using \eqref{eqn:c1c2}, we get 
\beqq{\tilde k = -\epsilon\Big(M\tilde m + M\tilde M M(AH^{-1} e+b))+o(\epsilon).}
Noting that $\kappa^* = -M(AH^{-1} e+b)$ leads to the desired claim in \eqref{eqn:lageps}. The proof of the lemma is hence complete.


\section{Algorithm for inequality constrained dynamic optimization and determination of optimal Lagrange multipliers}\label{app:MOM}

\begin{algorithm}
	\caption{Method of Multipliers-based Dynamic Optimization}
	\label{Alg_DP_aMoM}
	\begin{algorithmic}[1]
		\STATE Define $f_t(x_t,u_t)$ , $h_t(x_t,u_t)$ and $c_t(x_t,u_t)$ of the optimization problem for $t=0,1,\dots,T$.
		\FOR{$t=T$ to $0$}
		\FOR{every discretized value of $x_t \in \mathcal{X}$}
		\STATE Initialize Lagrange multiplier $\mu^0_t$ and penalty coefficient $c^0_p$. Define terminal cost $V_{T+1}(x)$. Set iteration $i=0$.
		\WHILE{$\|\mu^{i+1}_t-\mu^{i}_t\|\geq \epsilon$}
		\STATE Construct augmented Lagrangian function $L_{c^i_p}(x_t,u_t,\mu^i_t)$ (expression provided below).
		\STATE Solve $$u^i_t = \arg \min_{u_t\in\ALP U} \;\; L_{c^i_p}(x_t,u_t,\mu^i_t).$$
		\STATE Update Lagrange multiplier and penalty coefficient \beqq{\mu^{i+1}_t &= \max\{0, \mu^i_t+ c^i_p h_t(x_t,u^i_t) \}, \\ c^{i+1}_p &= \beta c^i_p, \beta \geq 1, \\ i &= i+1.}
		\ENDWHILE
		\STATE Store optimal control action $(u_t^i \to \gamma^*_t(x_t))$ and Lagrange multiplier $(\mu^i_t \to \mu^*_t(x))$. \cite{bertsekas2016nonlinear}
		\ENDFOR
		\STATE Store closed-loop optimal strategy $\gamma^*_t(x)$ and $\mu^*_t(x)$.
		\ENDFOR
	\end{algorithmic}
\end{algorithm}

The augmented Lagrangian is defined as
\beqq{& L_{c^i_p}(x_t,u_t,\mu^i_t)=  c_t(x_t,u_t) + V_{t+1}(f_t(x_t,u_t))+ \frac{1}{2c^i_p}\times\\& \sum_{k=1}^r\left(\left(\max\{0,\mu_{t,k}^i+c^i_p h_{t,k}(x_t,u_t)\}\right)^2-(\mu^i_{t,k})^2\right) .}

\end{document}

%% file: 2_Section_Problem_Formulation.tex
Consider a discrete-time dynamic control problem having the form:
\beq{\label{eqn:defn_state_dyn}
x_{t+1} = f_t(x_t,u_t), \quad t = 0,1,\dots,T.}
where $t$ denotes the discrete time instant, $x_t \in \ALP X \subset \Re^n$ is the state, $u_t \in \ALP U \subset \Re^m$ is the input or control, and $f_t$ is a function that describes the state dynamics. The control and the state are constrained, and the constraint function $h_t:\ALP X\times\ALP U\rightarrow\Re^r$ is of the form:
\beq{ \label{eqn:defn_constr}
 h_t(x_t,u_t)\leq 0 \quad \text{ for all } t= 0,1,\ldots,T.}
An admissible control map of the controller at time $t$ is a map $\gamma_t:\ALP X\rightarrow\ALP U$ such that $h(x,\gamma_t(x))\leq 0$ for all $x\in\ALP X$. We let $\Gamma_t$ denote the set of all such maps. The collection of admissible control maps denoted by $\VEC \gamma:=(\gamma_0,\ldots,\gamma_T)$, is referred to as the strategy of the controller. Let $\Gamma:=\prod_{t=0}^T\Gamma_t$ denote the set of all strategies of the controller. Note that a strategy is said to be open-loop if for all time $t\in\{0,\ldots,T\}$, there exists $u_t\in\ALP U$ such that $\gamma_t(x) = u_t$ for all $x\in\ALP X$.

The controller aims at minimizing a performance (cost) index:
\beq{ \label{eqn:defn_obj_fn}
J(\VEC \gamma) = c_{T+1}(x_{T+1}) + \sum_{t = 0}^{T} c_t(x_t,u_t),}
where $c_t:\ALP X\times\ALP U \rightarrow\Re$ is the per stage cost function. Let $\VEC\gamma^*$ denote the optimal strategy and $V_t$  denote the value function at time $t$ of the constrained dynamic optimization thus formulated. Along the optimal trajectory, let $\mu^*_t(x_t)$ denote the Lagrange multiplier corresponding to the inequality constraint and $\ALP A_t(x_t)$ denote the set of active constraints at $(x_t,\gamma^*_t(x_t))$, which is defined as
\beqq{\ALP A_t(x_t) = \Big\{j\in\{1,\ldots,r\}: \mu^*_{t,j}(x_t) >0\Big\}.}

The objective of this work is to define a solution for the perturbed optimization problem, where either the objective function or the constraints are perturbed w.r.t. the solution of the original problem \eqref{eqn:defn_state_dyn}-\eqref{eqn:defn_obj_fn}. The new optimization problem attempts at minimizing $J_\epsilon$, defined as:
\beqq{J_\epsilon(\VEC \gamma) = & c_{T+1}(x_{T+1}) + \sum_{t = 0}^{T} c_t(x_t,u_t) \\
& +\epsilon\left( \tilde c_{T+1}(x_{T+1}) + \sum_{t = 0}^{T} \tilde c_t(x_t,u_t) \right)\\
\text{ subject to } & h_t(x_t,u_t)+\epsilon \tilde h_t\leq 0,}
where $\tilde c_t$ and $\tilde h_t$ are perturbations of the objective function and the constraint function, respectively, and $\epsilon$ is used to parametrize the entire optimization problem for ease of analysis.

The following assumptions are introduced on the dynamic optimization problem.
\begin{assum}\label{assm:all}
For the unperturbed dynamic optimization problem,
\begin{enumerate}
    \item The cost function $c_t$, the constraint function $h_t$, and the state transition function $f_t$ are twice differentiable with continuous second derivative.
    \item The value function $V_t$ is twice differentiable with continuous second derivative for all $t$.
    \item For every time $t$, the one-step optimization problem (when using dynamic programming) satisfies regularity. In other words, the matrix $A$ defined as $A = [\nabla_{u_t} h_{t,i}(x_t,\gamma^*_t(x_t))]_{i\in\ALP A_t(x_t)}$ is a full rank matrix. 
    \item The matrix $\nabla^2_{u_tu_t} (c_t+V_{t+1}\circ f_t)(x_t,\gamma^*_t(x_t))$ is a positive definite matrix.  
\end{enumerate}
In addition, the perturbation $\epsilon$ is so small that the above conditions are satisfied by the perturbed dynamic optimization problem as well.\hfill $\Box$
\end{assum}
Note that these conditions are stronger than what is needed for the existence of a Lagrange multiplier. These assumptions can be relaxed by carefully constructing another equivalent cost function leading to an augmented Lagrangian method (see Section 4.2 of \cite{bertsekas2016nonlinear}), and would result in a cumbersome analysis. To ease exposition, such an analysis is not considered in this work.


%% file: 3_Section_Constraint_Perturbation_Warm_Start.tex
Consider a simple constrained optimization problem. Let $\ALP Z\subset\Re^p$ be a set and $g,\tilde g:\ALP Z\rightarrow \Re$ be cost functions. Let $q,\tilde q:\ALP Z\rightarrow\Re^s$ be the constraint functions. It is assumed that all these functions are twice differentiable with continuous and bounded second derivative in their domain $\ALP Z$. 

Consider the optimization problem:
\beq{\min_{z\in\ALP Z}\;\; & g(z) \nonumber \\ \text{ subject to } & q(z)\leq 0.\label{eqn:gopt}}
Let $z^*$ be the optimal solution and $\mu^*$ be the corresponding Lagrange multiplier. From the theorem of Lagrange multipliers, the following result holds:
\beqq{\nabla g(z^*) = -\nabla q(z^*)\mu^*.}
Now consider the following optimization problem:
\beq{\min_{z\in\ALP Z}\;\; & g(z)+\epsilon\tilde g(z) \nonumber \\ \text{ subject to } & q(z)+\epsilon\tilde q\leq 0.\label{eqn:tildegopt}}
Let $(\tilde z^*_\epsilon,\tilde \mu^*_\epsilon)$ denote the optimal solution and corresponding Lagrange multiplier pair and let $\tilde v^*_\epsilon$ denote the optimal value of the optimization problem.

\subsection{Solution of Perturbed Quadratic Programming}\label{sub:quadpro}
Consider the following quadratic optimization problem:
\beq{\label{eqn:gquadopt} \min_{z\in\Re^p}\;\; g(z) = & \frac{1}{2}z^\transpose H z+e^\transpose z \nonumber \\ \text{ subject to } & Az -b= 0,}
where $H$ is a positive definite matrix, and $e,A,b$ are appropriate constants. The solution to this optimization problem is given by:
\beq{\kappa^* & = -(AH^{-1}A^\transpose)^{-1} \Big(AH^{-1}e+b\Big) \label{eqn:quadprogkappa}\\
z^* &= -H^{-1}\Big(e+A^\transpose \kappa^*\Big).\label{eqn:quadprogzstar}}
The optimal value of the optimization problem is given by:
\beq{v^* := g(z^*) = \frac{1}{2}\left[\kappa^{*\transpose}(AH^{-1}A^\transpose) \kappa^* -e^\transpose H^{-1} e\right].\label{eqn:quadprogval}}
Let $\tilde g(z) = \frac{1}{2}z^\transpose \tilde H z+\tilde e^\transpose z$ and $\tilde q = -\tilde b$. Then, the new problem is formulated as:
\beqq{\min\;\; & \frac{1}{2}z^\transpose (H+\epsilon\tilde H) z+(e+\epsilon\tilde e)^\transpose z\\
\text{ subject to } & Az = b+\epsilon\tilde b. }
Naturally, this problem can be solved using the result in \eqref{eqn:quadprogkappa} and \eqref{eqn:quadprogzstar}. The aim however, is to identify an approximate solution to the perturbed quadratic optimization problem using the results of the original problem and the values of $\tilde H,\tilde e$ and $\tilde b$. To do so, the following lemma is needed.
\begin{lem}\label{lem:quadapprox1}
Let 
\begin{subequations}
\label{eqn:MMM}
\beq{& M = (AH^{-1}A^\transpose)^{-1},\quad \tilde M = AH^{-1}\tilde HH^{-1}A^\transpose,\\
& \tilde m = AH^{-1}\tilde e+\tilde b-AH^{-1}\tilde H H^{-1} e\\
& \tilde k = M \Big(\tilde M \kappa^*-\tilde m \Big)}
\end{subequations}
Then, the following expressions hold:
\beq{& (H+\epsilon\tilde H)^{-1} = H^{-1} -\epsilon H^{-1}\tilde H H^{-1} + o(\epsilon),\label{eqn:heps}\\
& \Big(A(H+\epsilon\tilde H)^{-1}A^\transpose\Big)^{-1} = M + \epsilon M\tilde MM + o(\epsilon),\label{eqn:ahaeps}\\
& \tilde\kappa^* =  \kappa^*+\epsilon\tilde k+o(\epsilon)\label{eqn:lageps}}
\end{lem}
\begin{pf}
See Appendix \ref{app:quadapprox1}. \hfill \QED
\end{pf}
Using the lemma above, the main result of this section is presented below. 
\begin{thm}\label{thm:quadprogperturbopt}
Define $d = -H^{-1}\Big(A^\transpose \tilde k+\nabla \tilde g(z^*)\Big)$. Then:
\beq{\tilde z^*_\epsilon - z^* = \epsilon d+ o(\epsilon).\label{eqn:peroptquadsol}}
Further, the change in value is: 
\beq{\tilde v^*_\epsilon - v^* = \epsilon\Bigg(& \tilde g(z^*)+\nabla g(z^*)^\transpose d\Bigg)+o(\epsilon) \label{eqn:peroptquadval}}
\end{thm}
\begin{pf}
The proof of the first statement is a direct application of \eqref{eqn:c1c2} (see Appendix \ref{app:quadapprox1}) to the expression in \eqref{eqn:quadprogzstar}. Indeed, the following must be noted first: \beqq{& -(H+\epsilon\tilde H)^{-1}\Big(e+\epsilon \tilde e+A^\transpose \tilde \kappa^*\Big)\\
& = -(H^{-1} -\epsilon H^{-1}\tilde H H^{-1} + o(\epsilon))\\
&\quad  \Big((e+A^\transpose \kappa^*)+\epsilon (\tilde e+A^\transpose \tilde k) +o(\epsilon)\Big).}
Now, \eqref{eqn:c1c2} and the expression for $z^*$ from \eqref{eqn:quadprogzstar} are used to arrive at the expression in \eqref{eqn:peroptquadsol}. The change in value \eqref{eqn:peroptquadval} is obtained as follows:
\beqq{\tilde v^*_\epsilon - v^* 
& = g(\tilde z^*_\epsilon)+\epsilon\tilde g(\tilde z^*_\epsilon) - g(z^*)\\
& = \epsilon\Big(\nabla g(z^*)^\transpose d+\tilde g(z^*)\Big)+o(\epsilon).}
The proof is hence complete.  \hfill \QED
\end{pf}
An important point to note in the expression for perturbations of the optimal value is that there is only one term that depends on the vector $d$. This fact will be key in extending the results obtained above to a perturbed DP problem.

\begin{exmp}
Consider a simple quadratic optimization problem with:
\beqq{H & = \begin{bmatrix} 6 & 2 \\ 2 & 1 \end{bmatrix}, \tilde H = \begin{bmatrix} 0.5 & 0.1 \\ 0.1 & 0.3 \end{bmatrix},
e = \begin{bmatrix} 3 \\5 \end{bmatrix},\tilde e = \begin{bmatrix} -0.4 \\ 0.2 \end{bmatrix}  \\
A & = \begin{bmatrix} 0 & -1\end{bmatrix},\quad  b = 0, \quad \tilde b = -0.3,\quad \epsilon = 1. }
The exact solutions to the optimization problem above are:
\beqq{z^* & = \begin{bmatrix} -0.5 \\0 \end{bmatrix},\qquad  \kappa^* = 4,\\
\tilde z^*_\epsilon & = \begin{bmatrix} -0.62 \\0.3 \end{bmatrix},\quad \tilde\kappa^* =  4.288.}
Using the expressions in \eqref{eqn:MMM} and in Theorem \ref{thm:quadprogperturbopt}, $\tilde k = 0.2$ and $d =  \begin{bmatrix} -0.125 \\0.3 \end{bmatrix}$ are computed. It is seen that the difference in the solutions $\tilde z^*_\epsilon - z^*$ is close to the approximation $d$ (recall $\epsilon = 1$ here). Further, it is observed that $\tilde\kappa^* -\kappa^*$ is also close to $\tilde k$ computed above.\hfill $\Box$
\end{exmp}

For the case where $z^* = 0$, there is a simplification in the expression for $d$, which is presented in the next corollary. This result is crucial for obtaining the expression for the WASP algorithm.

\begin{cor}\label{cor:Bw}
Define the matrix $B$ and vector $w$ as
\beq{B & = H^{-1}\begin{bmatrix} A^\transpose M A H^{-1}- I\Big| A^\transpose M  \end{bmatrix}, \quad w & = \begin{bmatrix} \tilde e\\\tilde b \end{bmatrix} \label{eqn:Bw}}
Let $d := Bw$. If $z^* = 0$, then the first order perturbation satisfies
\beqq{\tilde z^*_\epsilon &= \epsilon d+o(\epsilon)\\
\tilde \kappa^*_\epsilon &= \kappa^*+\epsilon (-M(AH^{-1}\tilde e + \tilde b))+o(\epsilon).}
\end{cor}
\begin{pf}
The proof follows from simple algebraic manipulations noting that the right side of \eqref{eqn:quadprogzstar} is equal to 0 when $z^* =0$. To see this, note that $-A^\transpose \kappa^* = e$. Using this expression in \eqref{eqn:lageps}, the following result is obtained:
\beqq{\tilde k &= M(-AH^{-1}\tilde H H^{-1}e - AH^{-1}\tilde e - \tilde b + AH^{-1}\tilde H H^{-1} e) \\
& = -M(AH^{-1}\tilde e + \tilde b).}
Further, since $z^* = 0$, $\nabla \tilde g(z^*) = \tilde e$. Theorem \ref{thm:quadprogperturbopt} is now used to get
\beqq{d = H^{-1}(A^\transpose M AH^{-1}\tilde e + A^\transpose M \tilde b) -H^{-1} \tilde e. }
This concludes the proof.  \hfill \QED
\end{pf}
An important observation from examining \eqref{eqn:Bw} is that if $z^* = 0$, then the first order perturbation in the optimal solution depends on the perturbation merely through $\tilde e$ and $\tilde b$. The matrix $B$ is dependent only on the original problem parameters for the quadratic program \eqref{eqn:gopt}, and not on the parameters of the perturbed problem.

Another special case is the situation when there is no equality constraint. For this case, the approximation is given below.
\begin{cor}\label{cor:quadprognocons}
Assume that there is no equality constraint. Then, 
\beqq{\tilde z_\epsilon^*-z^* &= -\epsilon H^{-1} \nabla \tilde g(z^*) + o(\epsilon),\\
\tilde v^*_\epsilon - v^* &= \epsilon\Big(\tilde g(z^*)-\nabla g(z^*)^\transpose H^{-1} \nabla \tilde g(z^*)\Big)+o(\epsilon).}
Further, if $z^* = 0$, then
\beq{\label{eqn:zs0nocons}\tilde z_\epsilon^*-z^* &= -\epsilon H^{-1} \tilde e+ o(\epsilon),\\
\label{eqn:zs0noconsval}\tilde v^*_\epsilon - v^* &= -\epsilon\Big(\tilde e^\transpose H^{-1} \tilde e\Big)+o(\epsilon).}
\end{cor}
\begin{pf}
 The proof is immediate.  \hfill \QED
\end{pf}

The first order perturbation in the optimal solution and the value function for various quadratic optimization problems have been established in this section. In the following subsection, these results are exploited to derive the same for a general nonlinear static optimization problem.

\subsection{Perturbation of General Nonlinear Optimization}
Consider the optimization problem in \eqref{eqn:gopt} with a nonlinear objective function $g$ and nonlinear inequality constraint function $q$. Assume that $g$ is twice differentiable with continuous second derivative and $q$ is once differentiable with continuous derivative.

Let $\ALP D = \{d\in\Re^p: z^*+d\in \ALP Z\}$ be the set of feasible directions at $z^*$. The optimization problem in \eqref{eqn:tildegopt} can be equivalently written as:
\beqq{\min_{d\in\ALP D} \;\;& g(z^*+d)+\epsilon\tilde g(z^*+d) \\
\text{ subject to } & q(z^*+d)+\epsilon\tilde q\leq 0.}
On taking a second order Taylor series expansion for the objective function,
\beqq{\hat L_\epsilon(d) =& g(z^*) +\nabla g(z^*)^\transpose d+\frac{1}{2}d^\transpose \nabla^2 g(z^*)d + \epsilon\tilde g(z^*)\\
& +\epsilon\nabla \tilde g(z^*)^\transpose d+\frac{\epsilon}{2}d^\transpose \nabla^2 \tilde g(z^*) d+o(\|d\|^2).}
Since $g(z^*)$ and $\tilde g(z^*)$ are constants and do not depend on $d$, these terms can be removed from the objective function. Further, the higher order terms are ignored to obtain a new objective function:
\beqq{L_\epsilon(d) = \frac{1}{2}d^\transpose \nabla^2 g(z^*)d+\frac{\epsilon}{2}d^\transpose \nabla^2 \tilde g(z^*) d\\
+\nabla g(z^*)^\transpose d +\epsilon\nabla \tilde g(z^*)^\transpose d}
The same operation is applied to the constraint, and results in the following optimization problem:
\beq{\min_{d\in\ALP D} \;\;& L_\epsilon(d)\nonumber\\
\text{ subject to } &  \nabla q(z^*)^\transpose d+(q(z^*)+\epsilon\tilde q)\leq 0.\label{eqn:lepsineq}}
This is still an inequality constrained quadratic optimization problem. Note that if some constraints are active (respectively inactive) and $d$ is small, those constraints can be expected to remain active (respectively inactive) at $z^*+d$. Let $\ALP A(z^*)$ denote the set of active constraints at $z^*$ at which the Lagrange multiplier is positive. Thus,
\beqq{\ALP A(z^*) &= \Big\{i\in\{1,\ldots,s\}: \mu_i^* >0\Big\}\\
& \subset \Big\{i\in\{1,\ldots,s\}: q_i(z^*) = 0\Big\}}
Such a hypothesis results in the following equality constrained quadratic optimization problem with equality constraints:
\beq{& \min_{d\in\ALP D}\;\; L_\epsilon(d)\nonumber\\
& \text{ s.t. } \nabla q_i(z^*)^\transpose d+(q_i(z^*)+\epsilon\tilde q_i) = 0 \text{ for all } i\in\ALP A(z^*).\label{eqn:Leps}}
This optimization problem can be solved using an iterative procedure if $\ALP D$ is a convex set using gradient projection method, conditional gradient method, or manifold suboptimization method (see Ch. 3 of \cite{bertsekas2016nonlinear} for a discussion on such algorithms). The problem when $\ALP D = \Re^p$ is studied here; this would be the case if $\ALP Z$ is an unconstrained set or $z^*$ is in the interior of the set $\ALP Z$.

Let $d^*_\epsilon$ be the optimal solution to the optimization problem above when $\ALP D = \Re^p$. When $\epsilon$ is small, $d^*_\epsilon$ is also expected to be small (indeed, for $\epsilon = 0$, $d^*_0 = 0$), and thus, the approximations are valid. Now, define:
\begin{subequations}
\label{eqn:staticoptmatrices}
\beq{H & = \nabla^2 g(z^*), && \tilde e = \nabla \tilde g(z^*),\\
A &= [\nabla q_i(z^*)^\transpose]_{i\in\ALP A(z^*)}, && \tilde b = -[\tilde q_i]_{i\in\ALP A(z^*)}.}
\end{subequations}
Using the results from Corollary \ref{cor:Bw}, $d^*_\epsilon$ can be computed approximately. In particular, the following observation is made.
\begin{thm}\label{thm:staticopt}
The optimal solution to the optimization problem in \eqref{eqn:Leps} is given by
\beqq{d_\epsilon^* = \epsilon Bw+o(\epsilon),}
where $B$ and $w$ are defined in \eqref{eqn:Bw} for $H,A,\tilde e,\tilde b$ as defined in \eqref{eqn:staticoptmatrices}. Further, the difference in optimal values are given by
\beqq{\tilde v^*_\epsilon - v^* = \epsilon \tilde g(z^*) -\epsilon\Big(\tilde e^\transpose H^{-1} \tilde e\Big)+o(\epsilon).}
\end{thm}
\begin{pf}
It is clear that for $\epsilon = 0$, $d^*_0 = 0$. This result is substituted in Corollary \ref{cor:Bw} to arrive at the conclusion.  \hfill \QED
\end{pf}

\begin{rem}\label{rem:rem1}
Recall that the optimization problem in \eqref{eqn:Leps} is an equality constrained problem, which is easy to solve using the approach adopted in Subsection \ref{sub:quadpro}. This method can lead to significant inaccuracies in cases where some inequality constraints that are not active in the original problem become active in the perturbed problem (this is particularly troublesome when the $\epsilon$ is large). This issue is alleviated by directly solving the inequality constrained quadratic programming problem in \eqref{eqn:lepsineq}. Note that this problem cannot be solved in closed form.
\end{rem}

For illustration, the approximation result is now applied to a resource allocation problem.

\subsection{Application to a Resource Allocation Problem}\label{sub:dynresonestep}
Consider the following optimization problem, parametrized by $x>0$:
\beqq{\min_{0\leq u\leq x} g(u) = -c\ln(u)-s\ln(x-u),}
where $c,s>0$ are constants. It can be readily established that the constraints are not active at the optimal solution, and that the optimal solution is given by $u^* = \frac{c}{c+s}x$, which lies in the open interval $(0,x)$. Further, the optimal value is $v^* = g(u^*) = \xi-(c+s)\ln(x)$, where $\xi$ is a constant dependent on $c$ and $s$ given by:
\beq{\label{eqn:xi}\xi = -c\ln\left(\frac{c}{(c+s)}\right)-s\ln\left(\frac{s}{(c+s)}\right).}

Consider the minimization of the following perturbed objective function:
\beqq{g(u)+\epsilon\tilde g(u) = -(c+\epsilon\tilde c)\ln(u)-(s+\epsilon\tilde s)\ln(x-u),}
where $\tilde c$ and $\tilde s$ are sufficiently small (not necessarily positive) so that $(c+\epsilon\tilde c)>0$ and $(s+\epsilon\tilde s)>0$.

Using the same expression as above, the optimal solution takes the form:
\beqq{\tilde u^*_\epsilon &= \frac{c+\epsilon\tilde c}{(c+s)+\epsilon(\tilde c+\tilde s)}x,\\
\tilde v^*_\epsilon &= g(\tilde u^*_\epsilon)+\epsilon\tilde g(\tilde u^*_\epsilon) = \xi_\epsilon- [(c+s)+\epsilon(\tilde c+\tilde s)]\ln(x),}
where $\xi_\epsilon$ can also be computed as in \eqref{eqn:xi} with appropriate modifications. Using \eqref{eqn:staticoptmatrices},
\beqq{H & = \nabla^2 g(u^*) = \frac{(c+s)^2}{x^2}\left(\frac{1}{c}+\frac{1}{s}\right),\\  
\tilde e &= \nabla \tilde g(z^*) = -\frac{(c+s)}{x}\left(\frac{\tilde c}{c}-\frac{\tilde s}{s}\right).}
Finally, \eqref{eqn:zs0nocons} is used to obtain the following:
\beq{\label{eqn:dynresonestepdstar}d^*_\epsilon = -\epsilon\frac{\tilde e}{H} = \epsilon\left(\tilde cs-\tilde sc\right)\frac{x}{(c+s)^2}.}
Pick $\epsilon = 1$, $c = 5, s = 10, \tilde c = -0.3, \tilde s = 0.4$, then:
\beqq{\tilde u_\epsilon^* - u^* = -0.0221x, \quad d^*_\epsilon = -0.0222x.}
In addition, the change in value is estimated to be:
\beqq{\tilde v^*_\epsilon - v^* &= \xi_\epsilon^* - \xi^*-\epsilon (\tilde c+\tilde s)\ln(x) \\
& = -0.1841 -(\tilde c+\tilde s)\ln(x).}
This result can be compared with the estimate based on Theorem \ref{thm:staticopt}, calculated below.
\beqq{\tilde v^*_\epsilon - v^* & \approx -\epsilon \frac{\tilde e^2}{H}-\epsilon (\tilde c+\tilde s)\ln(x) \\
&+ \epsilon \left(-\tilde c\ln\left(\frac{c}{(c+s)}\right)-\tilde s\ln\left(\frac{s}{(c+s)}\right)\right)\\
& = -0.2007- (\tilde c+\tilde s)\ln(x) .}
Once again, the estimate in error is in excellent agreement with the theoretical predictions.

%% file: 4_Section_WASP.tex
The aim now is to determine the first order perturbation in optimal policies and value functions for the perturbed constraint dynamic optimization problem. Suppose that the optimization problem \eqref{eqn:defn_state_dyn}-\eqref{eqn:defn_obj_fn} has been solved to get $(\gamma^*_t)_{t=1}^T$ and the corresponding optimal value functions $(V_t)_{t=1}^{T+1}$. The method of multipliers, described in Algorithm \ref{Alg_DP_aMoM} in Appendix \ref{app:MOM}, can be used to solve the stated dynamic optimization problem. A significant benefit of this algorithm is that it also returns the optimal Lagrange multiplier along with the optimal strategy and value function.

At time $t=T+1, V_{T+1}(x) = c_{T+1}(x)$. Based on the general formulation of DP as a backward induction problem \cite{bertsekas2012dynamic}, the following minimization problem is solved at time $T$:
\beqq{V_T(x_T) = \min_{u_T\in\ALP U}\quad  & c_T(x_T,u_T)+V_{T+1}(f(x_T,u_T))\\
 \text{ subject to } & h_T(x_T,u_T)\leq 0.}
Let $\mu^*_T(x_T)$ be the Lagrange multiplier corresponding to the inequality constraint and $\gamma^*_T(x_T)$ be the optimal strategy. In what follows, $V_T$ is assumed to be twice differentiable.

Now consider the perturbed problem:
\beqq{\min_{u_T\in\ALP U}\quad  & c_T(x_T,u_T)+V_{T+1}(f_T(x_T,u_T))\\
& +\epsilon (\tilde c_T(x_T,u_T)+\tilde V_{T+1}(f_T(x_T,u_T)))\\
\text{ subject to } & h_T(x_T,u_T)+\epsilon\tilde h_T\leq 0.}
Let $\bar\gamma^*_T(x_T)$ be the optimal strategy of the perturbed problem and $\bar V_T(x_T)$ be the optimal value function. For every $x_T$, the first order perturbation can be determined as a function of $x_T$ using the expressions in Theorem \ref{thm:staticopt}. Define $g,\tilde g$, and $q$ for a fixed $x_T$ as:
\beqq{g(u_T) &= c_T(x_T,u_T)+V_{T+1}\circ f_T(x_T,u_T),\\
\tilde g(u_T) &= \tilde c_T(x_T,u_T)+\tilde V_{T+1}\circ f_T(x_T,u_T),\\
q(u_T) &= h_T(x_T,u_T),\qquad \tilde q = \tilde h_T.}
Note that from the result of Theorem \ref{thm:staticopt}, term $V_{T+1}\circ f_T$ needs to be differentiated twice with respect to $u_T$. The following lemma provides the formula required to compute this derivative.
\begin{lem}
Let $\ALP Y\subset \Re^n$, $\ALP X\subset\Re^m$ and $\ALP U\subset\Re^p$. Let $V:\ALP Y\rightarrow\Re$ be a differentiable function and consider the composite function $V\circ f:\ALP X\times\ALP U\rightarrow\ALP Y$. On using the chain rule for multivariate functions,
\beqq{\nabla_{u} V\circ f(x,u) & = \nabla_{u} f|_{(x,u)} \nabla_{y} V|_{f(x,u)},\\
\nabla^2_{uu} V\circ f(x,u) & = \nabla_{u} f|_{(x,u)} \nabla^2_{yy} V|_{f(x,u)} \nabla_{u} f|_{(x,u)}^\transpose\\
& +\sum_{i=1}^n\nabla_{y_i} V|_{f(x,u)} \nabla^2_{uu} f_i|_{(x,u)}.}
\end{lem}
\begin{pf}
The proof is straightforward and results from repeated application of the chain rule reviewed in \eqref{eqn:chainrule}.  \hfill \QED
\end{pf}
Based on the stated assumptions, $\nabla^2_{uu}[ c_{T-1}+V_T\circ f_{T-1}]$ is positive definite at the optimal solution $(x_{T-1},\gamma^*_{T-1}(x_{T-1}))$. The backward induction algorithm is now applied, and Theorem \ref{thm:staticopt} is repeatedly used to derive the expression for first order perturbations at every time step. The following theorem captures the main result, and Algorithm \ref{Alg_WASP} outlines the proposed approach.
\begin{thm}\label{thm:dpperturb}
Define 
\begin{subequations}\label{eqn:dpperturb}
\beq{H(x) & = \nabla^2_{uu} c_t(x,\gamma_t^*(x))+\nabla^2_{uu} V_{t+1}\circ f_t(x,\gamma_t^*(x)),\\ 
\tilde e(x) &= \nabla_u \tilde c_t(x,\gamma_t^*(x))+\nabla_u \tilde V_{t+1}\circ f_t(x,\gamma_t^*(x)),\\
A(x) &  = [\nabla_u h_{t,i}(x,\gamma_t^*(x))^\transpose]_{i\in\ALP A(x,\gamma^*_t(x))}, \\
\tilde b(x) & = -[\tilde h_{t,i}]_{i\in\ALP A(x,\gamma^*_t(x))}.}
\end{subequations}
Construct the matrices $w_t(x)$ and $B_t(x)$ using the expressions in \eqref{eqn:Bw}. Define $d_t(x)$ and $\tilde V_t$ as:
\beq{d_t(x) =& B_t(x) w_t(x)\label{eqn:dteq}\\
\tilde V_t(x) =&  \tilde c_t(x,\gamma^*_t(x)) \nonumber\\
& + \nabla_u (c_t+V_{t+1}\circ f)^\transpose(x,\gamma^*_t(x)) d_t(x).\label{eqn:vteq}}
Then, the first order perturbations are given by:
\begin{subequations}\label{eqn:bareq}
	\beq{\bar \gamma^*_t(x) &= \gamma^*_t(x)+\epsilon d_t(x)+o(\epsilon),\\
		\bar V_t(x) &= V_t(x)+\epsilon\tilde V_t(x)+o(\epsilon).}
\end{subequations}
\end{thm}
\begin{pf}
The proof follows from a simple application of the principle of mathematical induction and Theorem \ref{thm:staticopt}.  \hfill \QED
\end{pf}

\begin{algorithm}
	\caption{WASP Algorithm} 
	\label{Alg_WASP} 
	\begin{algorithmic}[1]
		\STATE Define $f_t(x_t,u_t)$ , $h_t(x_t,u_t)$ and $c_t(x_t,u_t)$ of the unperturbed optimization problem for $t=0,1,\dots,T$.
		\STATE Perform constrained DP using Algorithm \ref{Alg_DP_aMoM} to simultaneously determine the optimal strategy-Lagrange multiplier pair $(\gamma^*_t(x),\mu^*_t(x))$, and optimal value function $V_t$.
		\STATE Determine $\tilde{c}_t(x_t,u_t)$ and $\tilde{h}_t$ in the perturbed optimization problem.
		\STATE Define terminal cost of perturbed problem $\tilde{V}_{T+1}(x)$.
		\FOR{$t=T$ to $0$}
		\FOR{every $x_t \in \mathcal{X}$}
		\STATE Using \eqref{eqn:dpperturb}, calculate $H(x_t), \tilde{e}(x_t)$. For $i\in\ALP A(x_t,\gamma^*_t(x_t))$, determine $A(x_t), \tilde{b}(x_t)$. \label{Alg_thm7_1}
		\STATE Construct $d_t(x_t)$ using \eqref{eqn:Bw} and \eqref{eqn:dteq}, or using \eqref{eqn:dtineq}. \label{Alg_thm7_2}
		\STATE Use \eqref{eqn:vteq} to estimate the perturbation in the value function $\tilde{V}_t(x_t)$.
		\ENDFOR
		\STATE Use \eqref{eqn:bareq} to evaluate first-order estimates of the strategy $\bar \gamma^*_t(x)$ and value function $\bar V_t(x)$.
		\ENDFOR
	\end{algorithmic}
\end{algorithm}

Often in many real-world problems, some constraints that are inactive at the optimal solution in a region of the state space become active in adjacent regions. For these problems, the above approach can lead to inaccuracies (as discussed earlier in Remark \ref{rem:rem1}). The accuracy of the proposed approach is controlled by explicitly incorporating the inequality constraints. To this end, define
\begin{subequations}\label{eqn:dpperturb2}
\beq{\tilde H(x) & = \nabla^2_{uu}\tilde c_t(x,\gamma_t^*(x))+\nabla^2_{uu} \tilde V_{t+1}\circ f_t(x,\gamma_t^*(x)), \\
e(x) & = \nabla_u c_t(x,\gamma_t^*(x))+\nabla_u V_{t+1}\circ f_t(x,\gamma_t^*(x)),\\
\hat A(x) &  = \nabla_u h_t(x,\gamma_t^*(x))^\transpose, \\
\hat b(x) & = -h_t^\transpose(x,\gamma_t^*(x))- \epsilon\tilde h_t.}
\end{subequations}
The following inequality constrained quadratic program can then be solved, and yields $d_t(x)$:
\beq{d_t(x) = \arg\min_{d\in\Re^m} & \frac{1}{2} d^\transpose (H(x)+\epsilon\tilde H(x)) d \nonumber\\
 & + d^\transpose (e(x)+\epsilon\tilde e(x))\nonumber\\
\text{subject to } & \hat A(x) d  \leq \hat b(x).\label{eqn:dtineq}}
The perturbation in value function can be computed using \eqref{eqn:vteq}. This approach leads to satisfactory results, as verified from the simulations in the following section. 


\section{Illustrative Examples}
In this section, the results from Theorem \ref{thm:dpperturb} (and \eqref{eqn:dpperturb2}) are applied to a dynamic resource allocation problem and a simplified velocity tracking problem.

\subsection{Dynamic Resource Allocation}
Consider the dynamic resource allocation problem, in which $f_t(x_t,u_t) = x_t - u_t$ and the total cost is given by:
\beqq{c_t(x_t,u_t) = -C_t\ln(u_t),}
where $C_1,\ldots,C_T$ are positive constants. The terminal cost is:
\beqq{c_{T+1}(x_{T+1}) = -C_{T+1}\ln(x_{T+1}).}
Define $S_t = \sum_{s=t}^{T+1}C_t$ for $t=0,\ldots,T$. The following results on the optimal solution and value function follows from the example in Subsection \ref{sub:dynresonestep}:
\beqq{\gamma_t^*(x_t) = \frac{C_t}{S_t}x_t,\qquad V_t^*(x_t) = \xi_t-S_t \ln(x_t).}
Now assume that the perturbations are of the form (pick $\epsilon = 1$ for simplicity):
\beqq{\tilde c_t(x_t,u_t) = -\tilde C_t\ln(u_t), \tilde c_{T+1}(x_{T+1}) = -\tilde C_{T+1}\ln(x_{T+1}).}
Define $\tilde S_t = \sum_{s=t}^{T+1}\tilde C_t$ for $t=0,\ldots,T$. From \eqref{eqn:dynresonestepdstar}, it can be concluded that:
\beqq{d_t^*(x_t) = \frac{\tilde C_t S_{t+1} - C_t \tilde S_{t+1}}{S_t^2}x_t.}
As shown through the numerical example in Subsection \ref{sub:dynresonestep}, this first order estimate of the difference in the optimal strategies between the perturbed and original problems is reasonably accurate.

\subsection{Simplified Velocity Tracking Problem}
Consider an optimization problem for simplified cruise control applications in vehicles. The goal in this problem is to determine the optimal desired acceleration profile to track a set reference velocity, given constraints on the feasible acceleration.

\subsubsection{Problem Formulation}
For simplicity, the vehicle considered in this example is modeled as a point mass. The only state variable is its velocity, $x_t = [v_t] \in \mathcal{X} \subset \mathbb{R}$ and the control action is the desired acceleration, $u_t = [a_t^{des}] \in \mathcal{U} \subset \mathbb{R}$. The kinematic equation that describes the discretized state dynamics is as follows:
\begin{equation}
\label{eqn:num_sim_state_eqn}
\begin{aligned}
v_{t+1} = v_t + a_t^{des}\Delta t, \quad t=0,1,\dots,T,
\end{aligned}
\end{equation}
where $t$ is the discrete time instant and $\Delta t$ is the time step. The performance criteria can be formulated by imposing costs on the deviation from the reference velocity, and the required control effort. This optimization problem is thus cast as a Linear Quadratic Regulator (LQR) problem (i.e. quadratic cost function with linear constraints). The total cost is given by:
\begin{equation}
\label{eqn:num_sim_cost}
\begin{aligned}
c_t(x_t,u_t) = w_p\left(v_t-v_{ref} \right)^2 + w_e \left(a_t^{des} \right)^2,
\end{aligned}
\end{equation}
where the parameters $w_p, w_e$ are respectively the performance weight that sets the cost for deviating from the reference state $v_{ref}$, and the control input weight which is interpreted as a measure of efficiency (often directly proportional to the desired propulsion energy).

The optimization problem is subject to the constraints posed by the state transition equation \eqref{eqn:num_sim_state_eqn}, and the following box constraints:
\beq{\label{eqn:num_sim_constraints}
a_t \in \left[a_t^{min}, a_t^{max} \right],}
where $\{a_t^{min}, a_t^{max}\}$ are the respective minimum and maximum vehicle acceleration limits for comfort.

To demonstrate the developed approach, the following quantities are perturbed: the reference velocity to be tracked, the control input constraints, as well as some of the performance objectives. The resulting policy and value function estimates are then evaluated by comparing them against the respective unperturbed and perturbed optimal solutions (that are determined using DP). While a simplified scenario is constructed here for validation, this approach is applicable and can be extended to several other LQR-based problems as well \cite{alam2015experimental,asadi2010predictive,rengarajan2020energy}.

\subsubsection{Perturbation of Original Cost Function and Constraints}
Pick $\epsilon = 1$ for simplicity and consider perturbations of the cost function and constraints having the form:
\begin{equation}
\label{eqn:num_sim_cost_pert}
\begin{aligned}
\tilde{c}_t(x_t,u_t) &= \tilde{w}_p\left(v_t-\tilde{v}_{ref} \right)^2 + \tilde{w}_e \left(a_t^{des} \right)^2, \\
a_t &\in \left[\tilde{a}_t^{min}, \tilde{a}_t^{max} \right].
\end{aligned}
\end{equation}
These perturbations could arise due to environmental disturbances (such as traffic, which affects the reference velocity to be tracked) and user preferences (which can affect the acceleration limits and weight parameters in the cost function).

\subsubsection{Simulation Results}
First, the original optimization problem is solved using Algorithm \ref{Alg_DP_aMoM}. The perturbed optimization problem is then solved using the WASP algorithm. Its performance is evaluated by comparison against the solutions obtained from using Algorithm \ref{Alg_DP_aMoM} (termed \textit{base} DP for use in this section) for both the unperturbed and the perturbed optimization cases. In this way, the bounds of performance improvement are set.

For $\{v_0, v_{ref}\}= \{10,12\}\si{m/s}$, $\{a_t^{min}, a_t^{max}\}= \{-2,2\}\si{m/{s^2}}$ and parameters $\{w_p,w_e\} = \{5,1\}$, the following results are obtained for the original (unperturbed) $5$-step optimization problem. The selected weights prioritize better tracking performance than efficiency of effort applied. Here, the solution converges to the reference velocity within three steps while satisfying all the imposed constraints.

\begin{figure}[!htbp]
	\centering
	\includegraphics[width=\columnwidth]{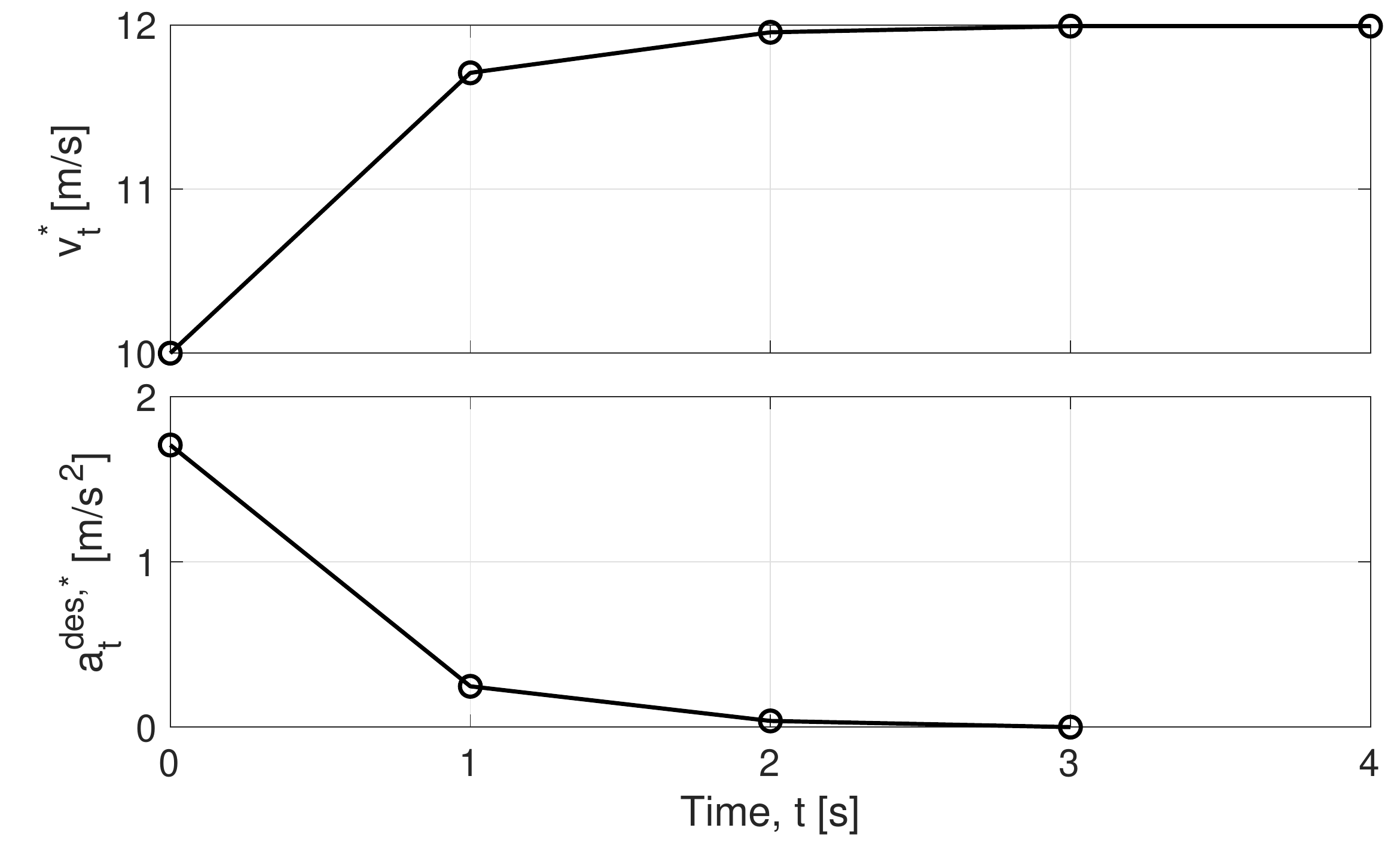}
	\caption{Optimal state and input trajectories for unperturbed velocity tracking problem using the \textit{base} DP algorithm.}
	\label{fig:num_sim_lqr_unpert_opt_st_inp}
\end{figure}

\begin{figure*}[!t]
	\centering
	\includegraphics[width=2\columnwidth]{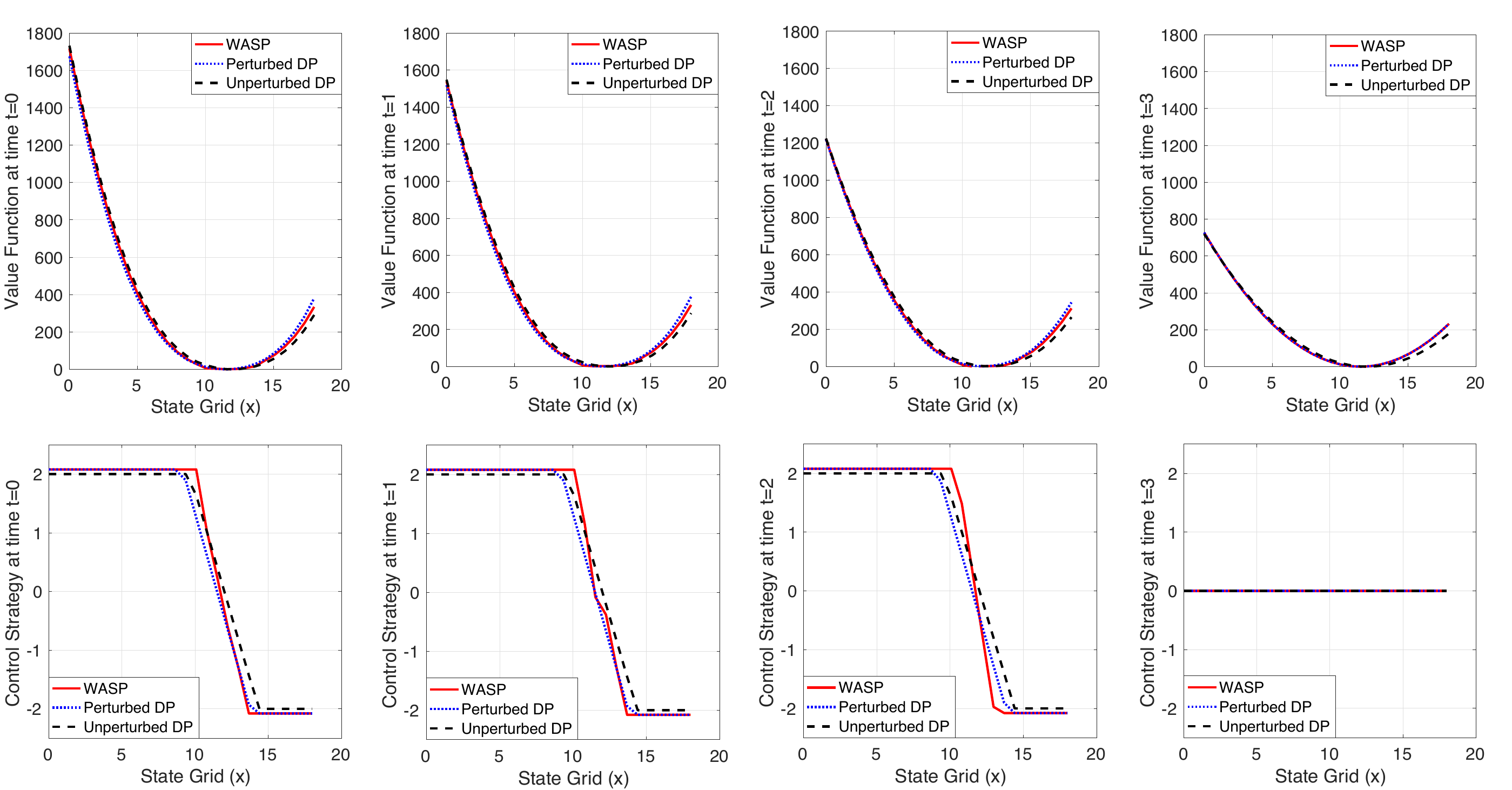}
	\caption{Value function and control strategy evolution for unperturbed and perturbed optimization cases demonstrating application of WASP algorithm for velocity tracking problem (standard deviation of injected Gaussian perturbation $=0.1$).}
	\label{fig:num_sim_lqr_pert_WASP_all_stddev_0.1}
\end{figure*}

Now consider \eqref{eqn:num_sim_cost_pert} in which constraints are perturbed with zero-mean Gaussian noise with standard deviation $0.1$ (in each of the upper/lower bounds for the action). Further, the reference velocity and weights in the cost function are perturbed randomly within $\pm \SI{10}{\%}$ of their original value. Figure \ref{fig:num_sim_lqr_pert_WASP_all_stddev_0.1} shows the comparison of the value function and control strategy trajectories from the unperturbed and the perturbed optimization cases. Note that in the simulations performed, \eqref{eqn:dpperturb2} and \eqref{eqn:dtineq} are used for computing the perturbed strategy and value function.

In the perturbed case, solutions from both the WASP and \textit{base} DP algorithms are visualized. It is seen that the estimates of the perturbed value function and control strategy from the WASP algorithm are in excellent agreement with the solution obtained using the \textit{base} DP for the perturbed optimization case. For $t=3$, the solutions from the WASP and \textit{base} DP for the perturbed optimization problem overlap with each other (as initialized). It is worth noting that the WASP algorithm does not use any of the results from the perturbed \textit{base} DP solution.

\begin{figure*}[!t]
	\centering
	\includegraphics[width=2\columnwidth]{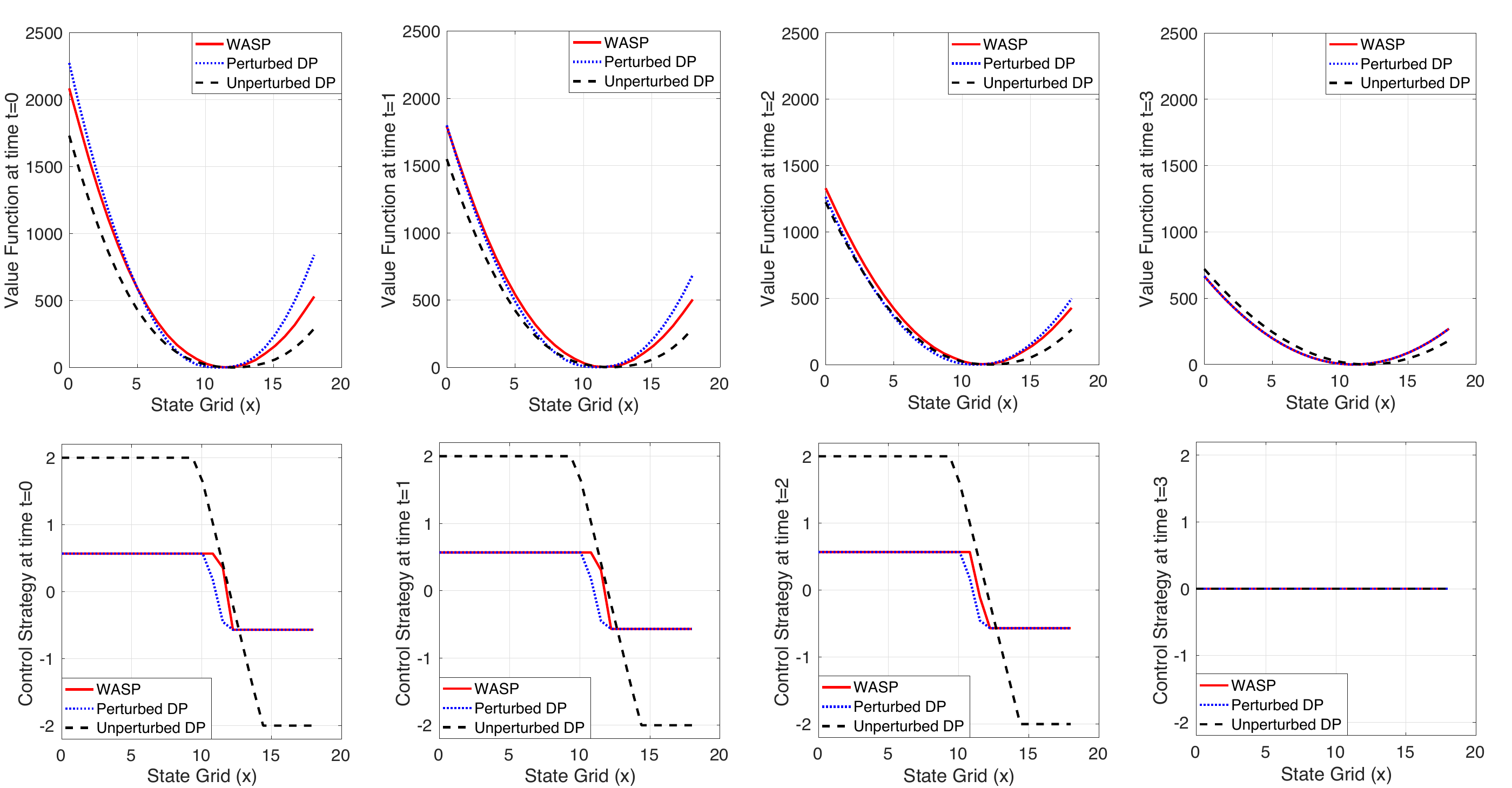}
	\caption{Value function and control strategy evolution for unperturbed and perturbed optimization cases demonstrating application of WASP algorithm for velocity tracking problem (standard deviation of injected Gaussian perturbation $=1.0$).}
	\label{fig:num_sim_lqr_pert_WASP_all_stddev_1.0}
\end{figure*}

The efficacy of the WASP algorithm can be clearly visualized when a much larger standard deviation is considered. Results from a simulation where the standard deviation is increased to $1.0$ are shown in Figure \ref{fig:num_sim_lqr_pert_WASP_all_stddev_1.0}. Here, the WASP algorithm always generates estimates of the value function that are closer to the true perturbed value than the unperturbed case. Further, the control strategy trajectories from $t=0$ closely follows the perturbed \textit{base} DP solution.


